\documentclass[a4paper,10pt, openany]{article}

\usepackage[french]{babel}
\usepackage[cp1252]{inputenc}
\usepackage{amsthm}
\usepackage{amsfonts}
\usepackage{amsmath}
\usepackage{amssymb}
\newcommand{\stackunder}[2]{\underset{#1}{#2}}

\newtheorem*{definition}{Définition}

\newtheorem{theorem}{Théorème}

\newtheorem{corollary}{Corollaire}

\begin{document}
\author{Stephane Collion}
\title{Transformation d'Abel et Formes Diff\'{e}rentielles Alg\'{e}briques }
\date{1996}
\maketitle

\begin{abstract}
Le but de cet article est de prouver qu'\'{e}tant donn\'{e}s un domaine
lin\'{e}airement concave $D$ dans l'espace projectif $\Bbb{CP}^{n}$, un
ensemble analytique complexe $V$ de dimension complexe 1 dans $D$ et une
1-forme m\'{e}romorphe $\varphi $ sur $V$, $V$ est inclus dans une courbe
alg\'{e}brique de $\Bbb{CP}^{n}$ et $\varphi $ est la restriction \`{a} $V$
d'une 1-forme alg\'{e}brique sur $\Bbb{CP}^{n}$ si et seulement si la
transform\'{e}e d'Abel $A(\varphi \wedge [V])$ du courant $\varphi \wedge
[V] $ est une 1-forme alg\'{e}brique sur l'espace dual $(\Bbb{CP}^{n})^{*}$,
o\`{u} une 1-forme alg\'{e}brique sur $\Bbb{CP}^{n}$ est une 1-forme
m\'{e}romorphe d\'{e}finie sur un rev\^{e}tement ramifi\'{e} de $\Bbb{CP}%
^{n} $. Ce r\'{e}sultat trouve son origine dans les th\'{e}or\`{e}mes d'Abel
inverse g\'{e}n\'{e}raux de Lie, Darboux, Saint-Donat, Griffiths et Henkin.

We prove in this article that given a linearly concave domain $D$ in the projective space $\Bbb{CP}^{n}$, a 1-dimensional comlex analytic set $V$ in $D$, and a meromorphic 1-form $\varphi$ on $V$, $V$ is a subset of an algebraic variety of $\Bbb{CP}^{n}$ and $\phi$ is the restriction to $V$ of an algebraic 1-form on $\Bbb{CP}^{n}$ if and only if the Abel transform $A(\varphi \wedge [V])$ of the analytic current $\varphi \wedge [V]$ is an algebraic 1-form on $(\Bbb{CP}^{n})*$, where an algebraic 1-form on $\Bbb{CP}^{n}$ is a meromorphic 1-form defined on a ramified analytic covering of $\Bbb{CP}^{n}$. This result has its origin in the general inverse Abel theorems of Lie, Darboux, Saint-Donat, Griffiths and Henkin.
\end{abstract}

\section{Rappel sur la transformation d'Abel. Notations}

\subsection{1-formes m\'{e}romorphes sur un ensemble analytique}

\smallskip Soit $V$ un ensemble analytique complexe de dimension complexe 1
dans un domaine $D\subset \Bbb{CP}^{n}.$ Une 1-forme $f$ sur $V$ est dite
m\'{e}romorphe (au sens de Griffiths) sur $V$ si il existe un ensemble
discret $S$ de points de $V$ avec $\mathit{Sing}V\subset S$ tel que $f$ est
une 1-forme holomorphe usuelle sur $V-S$ et tel que pour tout $w\in S$ il
existe un voisinage $U$ de $w$ dans $D$ et $m\in \Bbb{N}$ tel que:

\begin{center}
\[
\int\limits_{\left\{ z\in V\cap U:d(z,w)\geqslant r\right\} }f\wedge 
\overline{f}=O(r^{-m}) 
\]
\end{center}

Plus g\'{e}n\'{e}ralement, une 1-forme $\varphi $ \ sur $V$ sera dite
holomorphe sur $V$ (resp. m\'{e}romorphe sur $V$ \`{a} singularit\'{e}s dans
l'ensemble discret $S\subset V$, \textit{Sing}$V\subset S$ ), au sens de G.
Henkin ssi pour tout $z$ de $V$ il existe un voisinage $U$ de $z$ dans $D$
tel que:

$V\cap U=\left\{ z\in U:g(z)=0\right\} $ o\`{u} $g\in \mathcal{O}(U)$ et $dg$
est non identiquement nulle sur $V\cap U$, et : 
\[
\varphi _{_{\mid (V-S)\cap U}}=\text{\textit{Res}}_{_{\mid (V-S)\cap U}}%
\frac{f}{g} 
\]

o\`{u} $f\in \Omega ^{2}(U)$, (resp. $f\in \mathcal{M}^{2}(U)$ tel que
l'ensemble polaire $Z$ de $f$ v\'{e}rifie $Z\cap V\subset S$), o\`{u} 
\textit{Res} est le r\'{e}sidu de poincar\'{e}-Leray.

\subsection{Transform\'{e}e d'Abel}

\smallskip Un domaine $D$ de $\Bbb{CP}^{n}$ est dit\textit{\ l-}concave si
il peut s'\'{e}crire $D=\stackunder{\xi \in D^{*}}{\cup }\Bbb{CP}_{\xi
}^{n-1}$ o\`{u} $D^{*}$ est un domaine, dit dual de $D$, de l'espace dual $(%
\Bbb{CP}^{n})^{*}$ et o\`{u} l'on note $.\Bbb{CP}_{\xi }^{n-1}=\left\{ z\in 
\Bbb{CP}^{n}:\xi .z=\xi _{0}.z_{0}+\xi _{1}.z_{1}+...+\xi
_{n}.z_{n}=0\right\} $ pour $\xi \in D^{*}$. On utilisera ces notations pour
d\'{e}signer un domaine et son dual. Soit $V$ un ensemble analytique de
dimension complexe 1 dans $D$ \textit{l-}concave, et $\varphi $ une 1-forme
m\'{e}romorphe sur $V$. G\'{e}n\'{e}riquement, pour $\xi \in D^{*}$, $\Bbb{CP%
}_{\xi }^{n-1}$ coupe $V$ transversalement en un nombre fini de points $%
\left\{ Z^{1}(\xi ),...,Z^{N}(\xi )\right\} \subset V-S.$ Si en
coordonn\'{e}es affines $\varphi =\varphi _{1}dz_{1}+...+\varphi _{n}dz_{n}$%
, la transform\'{e}e d'Abel du produit ext\'{e}rieur $\varphi \wedge [V]$ de 
$\varphi $ avec le courant d'int\'{e}gration de Lelong $[V]$ est alors
d\'{e}finie pour $\xi \in D^{*}$ par: 
\[
A(\varphi \wedge [V])(\xi )=\sum\limits_{k=1}^{N}(Z^{k})^{*}\varphi (\xi
)=\sum\limits_{k=1}^{N}\left( \sum\limits_{j=1}^{n}\varphi _{j}(Z^{k}(\xi
))dZ_{j}^{k}(\xi )\right) 
\]

$A(\varphi \wedge [V])$ est une 1-forme m\'{e}romorphe dans $D^{*}$ ; elle
est holomorphe dans $D^{*}$ si $\varphi $ est holomorphe sur $V$.

\begin{theorem}
(Lie, Darboux, Saint-Donat et Griffiths).- Si $\varphi $ est non
identiquement nulle sur les composantes de $V$, alors $A(\varphi \wedge
[V])\equiv 0$ si et seulement si $V=\widetilde{V}\cap D$ o\`{u} $\widetilde{V%
\text{ }}$ est une courbe alg\'{e}brique de $\Bbb{CP}^{n}$ et $\varphi =%
\widetilde{\varphi }_{\mid V}$ o\`{u} $\widetilde{\varphi }$ est une 1-forme
rationnelle, holomorphe sur $\widetilde{V}$.
\end{theorem}

\begin{theorem}
(Henkin).- Si $\varphi $ est non identiquement nulle sur les composantes de $%
V$ et si $f=A(\varphi \wedge [V])\in \mathcal{M}^{1}(D^{*})$ se prolonge
\`{a} un domaine $\widetilde{D^{*}}\supset D^{*}$, alors il existe un
ensemble analytique $\widetilde{V}\subset \widetilde{D}$ et $\widetilde{%
\varphi }\in \mathcal{M}^{1}(\widetilde{V})$ tels que $V=\widetilde{V}\cap D$
et $\varphi =\widetilde{\varphi }_{\mid V}$.
\end{theorem}

Remarque: La d\'{e}monstration de ce th\'{e}or\`{e}me montre que l'ensemble
polaire de $f=A(\varphi \wedge [V])$ est de la forme $\bigcup\limits_{z\in
S}\left\{ \xi \in (\Bbb{CP}^{n})^{*}:\xi .z=0\right\} .$

\section{Formes Alg\'{e}briques, \'{e}nonc\'{e} du th\'{e}or\`{e}me
principal.}

Dans la th\'{e}orie classique, une fonction alg\'{e}brique sur $\Bbb{CP}^{n}$
est une restriction d'une fonction d\'{e}finie sur un rev\^{e}tement
ramifi\'{e} de $\Bbb{CP}^{n}$. On va d\'{e}finir de fa\c{c}on analogue les
formes diff\'{e}rentielles alg\'{e}briques.

Remarque 1: -Si $(X,p,Y)$, $p:X\rightarrow Y\,,\,$est un rev\^{e}tement
ramifi\'{e} d'une vari\'{e}t\'{e} complexe $Y$, un th\'{e}or\`{e}me de
Grauert et Remmert dit que $X$ est un ensemble analytique. Donc on peut
parler de formes m\'{e}romorphes sur $X$ au sens de Griffiths ou Henkin.

Remarque 2:-Soit $\mathcal{M}^{p}$ le faisceau des p-formes m\'{e}romorphes
sur une vari\'{e}t\'{e} complexe $Y\,$, $E\mathcal{M}^{p}$ l'espace
\'{e}tal\'{e} sur $Y$ associ\'{e}. Si $f_{a}\in \mathcal{M}_{a}^{p},f_{b}\in 
\mathcal{M}_{b}^{p}$ et si $\delta :I=[0,1]\rightarrow Y$ est un chemin de $%
Y $ de $a$ \`{a} $b$, $f_{b}$ est le prolongement de $f_{a}$ le long de $%
\delta $ si et seulement si il existe un rel\'{e}vement $\widehat{\delta }$
de $\delta $ \`{a} $E\mathcal{M}^{p}$ tel que $\widehat{\delta }(0)=f_{a}$
et $\widehat{\delta }(1)=f_{b}$. $E\mathcal{M}^{p}$ \'{e}tant
s\'{e}par\'{e}, on a l'unicit\'{e} d'un prolongement, ainsi que le
th\'{e}or\`{e}me de la monodromie.

\begin{theorem}
Soit $U$ un domaine de $\Bbb{CP}^{n}$ et $f\in \mathcal{M}^{1}(U).$ Les
trois propri\'{e}t\'{e}s suivantes sont \'{e}quivalentes:

\begin{enumerate}
\item  Il existe un rev\^{e}tement ramifi\'{e} $(X,p,\Bbb{CP}^{n})$
d'ensemble critique $Z$ et il existe $\widetilde{f}\in \mathcal{M}^{1}(X)$
au sens de Griffiths telle que $f=\sigma ^{*}\widetilde{f}$ pour une section 
$\sigma :U\rightarrow X$ de $p$.

\item  $f$ se prolonge le long de tout chemin $\delta $ de $\Bbb{CP}^{n}-Z$
o\`{u} $Z$ est un ensemble analytique de codimension complexe $\geqslant 1$,
de telle sorte qu'on obtienne qu'un nombre fini de prolongements en chaque
point et que $Z$ soit une singularit\'{e} alg\'{e}brique pour ces
prolongements (cf. ci-dessous).

\item  Si en coordonn\'{e}es affines $(z_{1},...,z_{n})$ dans $\Bbb{CP}^{n}=%
\Bbb{CP}_{\infty }^{n-1}+\Bbb{C}^{n}$, $f$ s'\'{e}crit $%
f=f_{1}dz_{1}+...+f_{n}dz_{n}$, alors $f_{1},...f_{n}$ sont alg\'{e}brique
au sur $\mathcal{M}(\Bbb{CP}^{n})$ au sens classique.
\end{enumerate}
\end{theorem}

Dans 2. $Z$ est une singularit\'{e} alg\'{e}brique si:

a/ $\forall m\in \Bbb{CP}^{n}-Z,\forall \widetilde{f}_{m}$ prolongement de$f$
en $m$, pour tout chemin ferm\'{e} $\delta $ de base $m$, $\exists l\in \Bbb{%
N}^{*}$ tel que le prolongement de $\widetilde{f}_{m}$ le long de $\delta
^{l}=\delta \vee \delta \vee ...\vee \delta $ (it\'{e}r\'{e} de $\delta l$
fois) redonne $\widetilde{f}_{m}$.

b/ $\forall m\in Z$, pour un voisinage $U$ de $m$ muni de coordonn\'{e}es $%
(z_{1},...,z_{n})$, les prolongements de $f$ dans $U-Z$ s'\'{e}crivant $%
f^{j}=f_{1}^{j}dz_{1}+...+f_{n}^{j}dz_{n}$ (dans une petite partie de $U-Z$%
), on a $\left| f_{k}^{j}(z)\right| =O(d(z,Z\cap U)^{-a_{jk}}$, o\`{u} $%
a_{jk}\in \Bbb{R}^{+}$.

\begin{definition}
- Soit $U$ un domaine de $\Bbb{CP}^{n}$ et $f\in \mathcal{M}^{1}(U)$. On dit
que $f$ est alg\'{e}brique sur $\Bbb{CP}^{n}$ si elle v\'{e}rifie l'une des
trois propri\'{e}t\'{e}s \'{e}quivalentes du th\'{e}or\`{e}me 3.
\end{definition}

Le principal r\'{e}sultat de cette note est alors:

\begin{theorem}
Soient donn\'{e}s au voisinage concave $U_{\xi _{0}}$ d'un hyperplan $\Bbb{CP%
}_{\xi _{0}}^{n-1}\,$de $\Bbb{CP}^{n}\,$un germe $V=V_{1}\cup ...\cup V_{k}$
d'ensemble analytique complexe de dimension 1 compos\'{e} de k composantes
irr\'{e}ductibles et $\varphi \in \mathcal{M}^{1}(V)$, $\varphi _{\mid
V_{j}}=\varphi _{j}$ non identiquement nulles. Soit $f=A(\varphi \wedge
[V])\in \mathcal{M}^{1}(U_{\xi _{0}}^{*})$. Alors les deux conditions
suivantes sont \'{e}quivalentes :

\begin{enumerate}
\item  $f$ est une 1-forme alg\'{e}brique sur $(\Bbb{CP}^{n})^{*}$ au sens
de la d\'{e}finition 4

\item  Il existe une courbe alg\'{e}brique $\widetilde{V\text{ }}\,$de $\Bbb{%
CP}^{n}$ telle que $V\subset \widetilde{V}$, et des 1-formes $\widetilde{%
\varphi }_{j}$ alg\'{e}briques sur $\Bbb{CP}^{n}$ telles que $\varphi _{j}=%
\widetilde{\varphi }_{j_{\mid V_{j}}}$ au sens de la d\'{e}finition 4.
\end{enumerate}
\end{theorem}

\textit{Remarque:- }Sous les hypoth\'{e}ses de 2., on montre que l'ensemble
de ramification (ou critique) de $f$ est de la forme $\bigcup\limits_{z\in
S}\left\{ \xi \in (\Bbb{CP}^{n})^{*}:\xi .z=0\right\} $ o\`{u} $S$ est un
ensemble dicret de point de $\Bbb{CP}^{n}$.

La d\'{e}monstration du th\'{e}or\`{e}me 5 sera donn\'{e}e en partie 4.
Remarquons d\'{e}ja qu'en consid\'{e}rant deux projection ind\'{e}pendantes
de $\Bbb{CP}^{n}-\Bbb{CP}^{n-3}$ sur deux plans $\Bbb{CP}^{2}$ puis en
consid\'{e}rant l'intersection des images inverses des prolongements obtenus
dans ces plans, il suffit de faire la d\'{e}monstration dans le cas n=2. Le
th\'{e}or\`{e}me 3 sera \'{e}galement d\'{e}montr\'{e} dans $\Bbb{CP}^{2}$
ce qui ne change rien.

\textit{Application:}

\begin{corollary}
-Soient $f^{1}(\xi _{0},\xi _{1})$,..,$f^{N}(\xi _{0},\xi _{1})$ N fonctions
m\'{e}romorphes dans un domaine $D^{*}$ de $\Bbb{CP}^{2}$ v\'{e}rifiant
chacune l'\'{e}quation diff\'{e}rentielle des ''ondes de choc'': 
\[
f^{j}.\frac{\partial f^{j}}{\partial \xi _{0}}=\frac{\partial f^{j}}{%
\partial \xi _{1}} 
\]
et soit $F=\sum\limits_{j=1}^{N}f^{j}$. Supposons que $dF$ soit
alg\'{e}brique, ce qui est le cas en particulier si $F$ est alg\'{e}brique.
Alors chaque $f^{j}$ est alg\'{e}brique.
\end{corollary}

\textit{D\'{e}monstration du corollaire:}- D'apr\'{e}s la proposition 3.2
dans l'article de G. Henkin, et en identifiant $D^{*}\,$au domaine dual d'un
domaine $D$ de $\Bbb{CP}^{2}\,$, on peut construire dans $D$ un ensemble
analytique $V\,$ tel que $\left\{ \left( 1:f^{j}(\xi ):z_{2}^{j}(\xi
)\right) \right\} _{j=1,...,N}=V\cap \Bbb{CP}_{\xi }^{1},\xi \in D^{*}$ et $%
dF=A(dz_{1}\wedge [V])$. Alors, d'apr\'{e}s le th\'{e}or\`{e}me 5, $V$ est
alg\'{e}brique, soit $V=\left\{ P(z_{1},z_{2})=0\right\} $ pour un
polyn\^{o}me $P$, donc $P(f^{j}(\xi ),-\xi _{0}-\xi _{1}.z_{1}^{j}(\xi
))=P(f^{j}(\xi ),-\xi _{0}-\xi _{1}.f^{j}(\xi ))=0$, d'o\`{u} $f^{j}$ est
alg\'{e}brique.

\smallskip

\textit{Remarque 1: }Reprenons les notations des th\'{e}or\`{e}mes 4 et 5 et
supprimons l'hypoth\'{e}se de croissance polynomiale de\textit{\ }$%
f=A(\varphi \wedge \left[ V\right] )$ vers $Z$: dans le th\'{e}or\`{e}me 4.1
on suppose seulement que $\widetilde{f}\in \mathcal{M}^{1}(X-p^{-1}(Z))\,$et
dans 2. on supprime la condition b/ (3. est supprim\'{e}e). La
d\'{e}monstration du th\'{e}or\`{e}me 5 montre que l'on a encore l'existence
d'une courbe alg\'{e}brique $\widetilde{V}$ contenant $V$. Ceci permet aussi
de renforcer le corollaire.

\textit{Remarque 2: }Il\textit{\ }faut comprendre le th\'{e}or\`{e}me 5
comme une caract\'{e}risation locale d'une courbe alg\'{e}brique. En
particulier, si l'on consid\`{e}re sur $V$ une 1-forme ''coordonn\'{e}e '' $%
dz_{j}$, la transform\'{e}e d'Abel $f^{j}=A(dz_{j}\wedge [V])\,$est la
diff\'{e}rentielle de la fonction $A_{V}^{j}(\xi
)=\sum_{k=1}^{N}Z_{j}^{k}(\xi )$ o\`{u} l'on a repris les notations du
paragraphe 1.2. Cette fonction $A_{V}^{j}(\xi )\,$est d\'{e}ja
\'{e}tudi\'{e}e par Blaschke, Bol, Saint-Donat, Wood et G.Henkin (elle est
appel\'{e}e par ce dernier transform\'{e}e d'Abel de $V$). Le r\'{e}sultat
obtenu par Henkin est que $V$ est la r\'{e}striction d'une courbe
alg\'{e}brique au domaine $U$ si et seulement si $A_{V}^{j}(\xi )$ est
rationnelle pour un j quelconque. Nous obtenons:

\begin{corollary}
- $V$ est inclus dans une courbe alg\'{e}brique si et seulement si $%
A_{V}^{j}(\xi )$ est alg\'{e}brique pour un $j$ quelconque (ou ce qui
revient au m\^{e}me par le corollaire 6, si et seulement si $f^{j}$ est alg%
\'{e}brique pour un $j$ quelconque).
\end{corollary}

\end{document}